\documentclass[11pt]{article}

\usepackage[margin=1in]{geometry}
\usepackage{amsmath,amssymb,amsthm,mathtools}
\usepackage{microtype}
\usepackage{enumitem}
\usepackage[colorlinks=true,linkcolor=blue,citecolor=blue,urlcolor=blue]{hyperref}
\usepackage{xcolor}

\newtheorem{theorem}{Theorem}[section]
\newtheorem{lemma}[theorem]{Lemma}

\newtheorem{proposition}[theorem]{Proposition}
\newtheorem{corollary}[theorem]{Corollary}

\newtheorem{conjecture}[theorem]{Conjecture}

\newcommand{\psr}{\operatorname{psr}}
\newcommand{\sat}{\operatorname{sat}}
\newcommand{\cG}{\mathcal{G}}

\title{The Plane-Saturation Ratio of Planar Graphs\\
	}
\author{ Qing Cui\footnote{School of Mathematics,
		Nanjing University of Aeronautics and Astronautics,
		Nanjing 210016,  China. {\tt cui@nuaa.edu.cn}},\ \ Shengjin Ji\footnote{School of Mathematics and Statistics, Shandong University of Technology,
		Zibo 255000, China. {\tt
			jishengjinjin@sdut.edu.cn}} \footnote{Corresponding author: Shengjin Ji}}
\date{}

\begin{document}
	\maketitle
	
	\begin{abstract}
		%Clifton and Salia introduced the plane-saturation ratio in order to measure how far a plane subgraph of a planar graph may deviate from an optimal embedding while still being saturated. 
%For a planar graph \(G\), a plane subgraph \(H\) is plane-saturated with respect to \(G\) if adding any edge to
%	the fixed drawing of \(H\), possibly with new vertices, either introduces a
%	crossing or causes the resulting graph  to no longer be a subgraph of
%	\(G\). The plane-saturation ratio \(\psr(G)\) is the minimum value of
%	\(e(H)/e(G)\) over all plane‐saturated subgraphs \(H\) of $G$. For nonnegative integers
%	\(k_1,k_2\), let \(\cG_{k_1,k_2}\) be the class of planar graphs in which, for
%	\(i\in\{1,2\}\), at most \(k_i\) vertices of degree \(i\) have any prescribed
%	neighborhood; the convention \(k_i=0\) forbids vertices of degree \(i\).
%	Clifton and Salia conjectured that every \(G\in\cG_{k_1,k_2}\) with at least
%	one edge satisfies
%	$\psr(G)>\frac{1}{9+k_1+6k_2}.$
%	We prove this conjecture. The proof first exploits the unlabeled nature of plane saturation to obtain structural restrictions, and then combines these restrictions with planar counting and the low-degree multiplicity bounds.

	%We study the plane-saturation ratio recently introduced by Clifton and Salia. 
	For a planar graph \(G\), a plane subgraph \(H\) is
	plane-saturated with respect to \(G\) if adding any edge to the fixed
	drawing of \(H\), possibly together with new vertices, either creates a
	crossing or results in a graph that is not a subgraph of \(G\). The
	plane-saturation ratio \(\psr(G)\) is the minimum of \(e(H)/e(G)\) over
	all plane-saturated subgraphs \(H\) of \(G\). For nonnegative integers
	\(k_1\) and \(k_2\), let \(\cG_{k_1,k_2}\) denote the class of planar
	graphs in which, for each \(i\in\{1,2\}\), at most \(k_i\) vertices of
	degree \(i\) have any prescribed neighborhood, where \(k_i=0\) means
	that vertices of degree \(i\) are forbidden. We prove that every
	\(G\in\cG_{k_1,k_2}\) with at least one edge satisfies
	$
	\psr(G)>\frac{1}{9+k_1+6k_2}.
	$
	This settles two conjectures of Clifton and Salia.

	\end{abstract}
	
	
	
	\section{Introduction}
The graph saturation problem asks for the minimum number of edges in an
\(n\)-vertex graph that contains no copy of a fixed graph \(F\), but in which
adding any missing edge creates a copy of \(F\). This minimum is the
saturation number \(\sat(n,F)\). The study of saturation was initiated by
Zykov \cite{Zykov}; Erd\H{o}s, Hajnal, and Moon determined the exact value of
\(\sat(n,K_t)\) \cite{ErdosHajnalMoon}, and Bollob\'as extended the complete
graph result to hypergraphs \cite{Bollobas}. Further general bounds and
related developments were obtained in
\cite{CameronPuleo,FaudreeGould,KaszonyiTuza}; see also the survey
\cite{FaudreeFaudreeSchmitt}.

Clifton and Salia \cite{CliftonSalia} introduced a variant of saturation for
partial plane embeddings of planar graphs. This direction is also related to
saturated topological drawings and to extension problems for partial planar
drawings; see, for example,
\cite{%AngeliniEtAl,
	ChaplickKluteParadaRollinUeckerdt,JelinekKratochvilRutter,KynclPachRadoicicToth}.
Related questions for maximal planar host graphs were studied in
\cite{BaratBlazsikKeszeghZheng,CliftonSimon}. The present paper continues
the study of Clifton and Salia and proves the full form of one of their main
conjectures.

Let \(G\) be a planar graph with at least one edge, and let \(H\) be a plane
graph that is isomorphic to a subgraph of \(G\). Following
\cite{CliftonSalia}, we say that \(H\) is \emph{plane-saturated with respect
	to \(G\)} if adding any edge to the fixed drawing of \(H\), possibly together
with new vertices, either creates a crossing or produces a graph that is no
longer isomorphic to a subgraph of \(G\). The \emph{plane-saturation ratio}
of \(G\) is
\[
\psr(G)=\min_H\frac{e(H)}{e(G)},
\]
where the minimum is over all plane-saturated subgraphs \(H\) of \(G\).
This is the unlabeled version of the problem: containment is considered up
to graph isomorphism, rather than with fixed vertex labels.

Clifton and Salia \cite{CliftonSalia} showed that \(\psr(G)\)
can be arbitrarily small for planar graphs \(G\).
Their examples contain arbitrarily many vertices of degree
one or two with the same neighborhood.
This led them to consider the following classes.
For nonnegative integers \(k_1,k_2\), let
\(\cG_{k_1,k_2}\) denote the class of planar graphs \(G\)
such that
\[
\bigl|\{v\in V(G):d_G(v)=i,\ N_G(v)=X\}\bigr|\le k_i
\]
for every \(i\in\{1,2\}\) and every \(X\subseteq V(G)\).
In particular, \(k_i=0\) means that \(G\) has no vertices
of degree \(i\).

Clifton and Salia proved the bound below for \(k_1\ge 1\)
and \(k_2\ge 0\), except when
\((k_1,k_2)\in\{(1,0),(2,0)\}\)
\cite[Theorem~1.5]{CliftonSalia}.
They conjectured that the same bound holds for all
nonnegative integers \(k_1,k_2\)
\cite[Conjecture~4.1]{CliftonSalia}. We restate it as follows. 

\begin{conjecture}\label{conj:CliftonSalia}
	For all nonnegative integers \(k_1,k_2\), every graph
	\(G\in\cG_{k_1,k_2}\) with at least one edge satisfies
	\[
	\psr(G)>\frac{1}{9+k_1+6k_2}.
	\]
\end{conjecture}

Our main result confirms this conjecture.

\begin{theorem}\label{thm:main}
	For all nonnegative integers \(k_1,k_2\), every graph
	\(G\in\cG_{k_1,k_2}\) with at least one edge satisfies
	\[
	\psr(G)>\frac{1}{9+k_1+6k_2}.
	\]
\end{theorem}

Taking \(k_1=k_2=0\), we obtain \(\psr(G)>1/9\) for
every planar graph \(G\) of minimum degree at least three.
Thus Theorem~\ref{thm:main} also settles
\cite[Conjecture~4.2]{CliftonSalia}.
For each fixed pair \((k_1,k_2)\), the bound is
asymptotically sharp by
\cite[Construction~2]{CliftonSalia}.

Our proof builds on the ideas of Clifton and Salia
\cite{CliftonSalia}. We analyze the tree components of
a spanning plane-saturated subgraph and combine the
resulting structural restrictions with the edge bound
for bipartite planar graphs and the conditions defining
\(\cG_{k_1,k_2}\) to obtain the desired lower bound.

	\section{Preliminaries}

		We introduce some notions necessary for the rest of the paper. For a graph \(G\), the set of its vertices is
	denoted by \(V(G)\) and the set of its edges by \(E(G)\). The numbers of
	vertices and edges are denoted by \(v(G)\) and \(e(G)\), respectively. Let $c(G)$ be components of $G$.  If
	\(G\) is a planar graph, then \(f(G)\) denotes the number of faces.  For a
	vertex \(v\) of \(G\), the neighborhood of \(v\) is the set of all vertices
	\(u\) such that \(vu\) is an edge of \(G\) and is denoted by \(N_G(v)\).
	The degree of \(v\) is the size of \(N_G(v)\) and is denoted by \(d_G(v)\).
	%When the host graph is apparent from the context, we may use \(N(v)\) and \(d(v)\) instead of \(N_G(v)\) and \(d_G(v)\), respectively. 
	When \(G\) is
	a graph and \(e\) is a pair of its nonadjacent vertices, \(G+e\) denotes
	the graph obtained from \(G\) by adding the edge \(e\). When \(U\) is a
	subset of vertices of \(G\), \(G[U]\) denotes the subgraph of \(G\) induced
	by \(U\), and we denote \(G[V(G)\setminus U]\) by \(G-U\). The cycle on
	\(n\) vertices is denoted by \(C_n\).  For
	\(V_0\subseteq V(G)\) and \(x\in V(G)\), put $N_{V_0}(x)=N_G(x)\cap {V_0}, d_{V_0}(x)=|N_{V_0}(x)|$. 
	A vertex \(x\in V(G)\) is called \(V_0\)-positive if
		\(d_{V_0}(x)>0\), and \(V_0\)-zero otherwise. %后面再检查那里还用了。都更正过来。%We write \((t)_+=\max\{t,0\}\). 
		Let $c(G)$ be components of $G$.   Then Euler's formula have the general version: $v(G)-e(G)+f(G)=1+c(G)$. In particular, $v(G)-e(G)+f(G)=2$ when $G$ is connected.

We present the following standard consequences of Euler's formula;
see~\cite[Chapter~4]{Diestel2017}.

\begin{proposition}\label{prop:planar-bounds}
	Let \(G\) be a simple planar graph with \(v(G)\ge3\). Then
	\[
	e(G)\le3v(G)-6.
	\]
	If \(G\) is bipartite, then
	\[
	e(G)\le2v(G)-4.
	\]
\end{proposition}

Applying Proposition~\ref{prop:planar-bounds} to the components of
order at least three and checking the remaining components directly
gives the following uniform bound.

\begin{corollary}\label{cor:component-planar}
	Every simple planar graph \(G\), including the empty graph, satisfies
	\begin{equation}\label{eq:component-planar}
		e(G)\le3\bigl(v(G)-c(G)\bigr).
	\end{equation}
	In particular, \(e(G)\le3v(G)-3\) whenever \(G\) is nonempty.
\end{corollary}

	 %可能还要加点别的符号
	%Observe that removing isolated vertices from the host graph \red{does not change \(\psr(G)\)}or  {\color{red}the degree-one and degree-two neighborhood multiplicities}, so we may also assume throughout the 	paper that $G$ has no isolated vertices. \cite{CliftonSalia}. %这句话可能要调整，degree-one表述不规范 Hence we assume that \(G\) has no isolated vertices.最后一句话有问题，删除孤立点可能会改变母图的psr值。这句话放这不合。
	We begin with some additional terminology and notation. 	
	Let \(G\)
	be a planar graph with $e(G)>0$ and \(H\) a plane-saturated subgraph of \(G\). Let \(I\)
	denote the set of isolated vertices of \(H\). We refer to the graph \(H-I\)
	as the \emph{skeleton}, denoted by $S$ and consequently refer to the vertices in
	\(V(H)\setminus I\) as \emph{skeleton vertices}. In particular, the
	skeleton is a planar graph without a specified embedding. %We call a skeleton component with exactly two vertices a \emph{matching component}.  The \emph{skeleton} \(S\) is the plane graph obtained from \(H\) by deleting the vertices of \(I\), with the inherited embedding of $H$. 

	We
	will sometimes refer to an \emph{embedding}
	\(\phi\colon V(H)\to V(G)\) of \(H\) in \(G\), where \(\phi\) is a
	bijective function such that, for every edge \(uv\) of \(H\),
	\(\phi(u)\phi(v)\) is an edge of \(G\). With a slight abuse of notation,
	for every edge \(uv\) of \(H\), we denote
	\(\phi(uv):=\phi(u)\phi(v)\). 
		%A\emph{bijective homomorphism} from \(H\) to \(G\) is a bijection \(\phi:V(H)\to V(G)\) such that $xy\in E(H) \Longrightarrow \phi(x)\phi(y)\in E(G)$.
		%{\color{red}Nonedges need not be preserved. Since plane saturation is unlabeled,if two distinct nonadjacent vertices $u,v$ are cofacial in the fixeddrawing of $H$ and $H+uv$ admits an embedding into $G$,	then $H$ is not plane-saturated.}
	
	  Recall one assertion  \cite{CliftonSalia} that adding isolated vertices to $H$ does not affect whether 
$H$ is a plane‐saturated subgraph of 
$G$ or not.	Hence, we may assume that $H$ is spanning.
%Let $	I=\{x\in V(H):d_H(x)=0\}.$ 

		%\red{	It suffices to minimize the saturation ratio over spanning drawings. Given a plane-saturated drawing, insert its missing vertices as isolated points and extend its edge-preserving map into $G$ over them.  If an edge could be added to the enlarged drawing, deleting the inserted vertices other than its endpoints would give	a permitted edge addition to the original drawing, possibly with new endpoints.  Thus the enlarged drawing remains saturated and has	the same edge count.}
		
		Observe that the assertion in \cite{CliftonSalia} that removing isolated
		vertices from $G$ does not affect $\psr(G)$ does not hold in general. For example, let $G=K_3\cup K_2\cup K_1$, and let $z$ be its isolated
		vertex. Then
		$
		\psr(G-z)=\frac{3}{4}<1=\psr(G).
		$
		 Fortunately, the following reduction is valid.
	\begin{lemma}[]\label{lem:reductions}
		If $Z$ is the set of isolated vertices of $G$, then
$
		\psr(G-Z)\le\psr(G).
$
	\end{lemma}
	
	\begin{proof}
		Let \(H\) be a spanning plane-saturated drawing with respect to \(G\), and
		fix an edge-preserving bijection
	$\phi:V(H)\to V(G).
		$
		Put
		$
		Z_H=\phi^{-1}(Z).
		$
		Since every vertex of \(Z\) is isolated in \(G\), every vertex of \(Z_H\)
		is isolated in \(H\). Let $
		H_0=H-Z_H.
		$
		Then \(H_0\) is a spanning drawing with respect to \(G-Z\), and
		$
		e(H_0)=e(H).
		$

		We claim that \(H_0\) is plane-saturated with respect to \(G-Z\). Suppose
		otherwise. Since \(H_0\) is spanning with respect to \(G-Z\), an admissible
		edge addition cannot introduce new vertices. Thus there are two nonadjacent
		vertices \(u,v\in V(H_0)\), lying on a common face of the drawing of \(H_0\),
		such that the edge \(uv\) can be drawn without crossings and \(H_0+uv\)
		admits an edge-preserving bijection into \(G-Z\). Since the vertices of
		\(Z_H\) are isolated in the drawing of \(H\), the arc representing \(uv\)
		may be chosen to avoid them. Extending the edge-preserving bijection
		\(H_0+uv\to G-Z\) by any bijection from \(Z_H\) onto \(Z\), we obtain an
		edge-preserving bijection from \(H+uv\) into \(G\),  contradicting the
		plane-saturation of \(H\). 
		
		Since \(e(G-Z)=e(G)\), we obtain
		\[
		\psr(G-Z)
		\le
		\frac{e(H_0)}{e(G-Z)}
		=
		\frac{e(H)}{e(G)}.
		\]
		Taking the minimum over all spanning plane-saturated drawings \(H\) with
		respect to \(G\) gives
	$
		\psr(G-Z)\le \psr(G).
		$
	\end{proof}

The example above shows that the inequality in
Lemma~\ref{lem:reductions} can be strict.

Now let $G\in\cG_{k_1,k_2}$, and let $Z$ be its set of isolated
vertices. Deleting $Z$ changes neither the degrees nor the
neighborhoods of the remaining vertices, so
$G-Z\in\cG_{k_1,k_2}$. By Lemma~\ref{lem:reductions}, proving the
lower bound in Theorem~\ref{thm:main} for $G-Z$ suffices to prove
it for $G$. Thus, in the proof of this lower bound, we may assume
that $G$ has no isolated vertices.

%In fact, by the example above, the inequality in Lemma~\ref{lem:reductions} can be strict. Let $G\in \cG_{k_1,k_2}$. Let $Z$ be the set of insolated vertices of $G$. By Lemma~\ref{lem:reductions}, we know that $\psr(G-Z)\le \psr(G)$.  to verify Theorem \ref{thm:main} can be reduced to verify that  the lower bound holds for $G-Z$. We may also assume throughout the  paper that $G$ has no isolated vertices.

Put
$
n=v(G)=v(H), $ $h=e(H), s=v(S),
\Lambda=9+k_1+6k_2.$
 Note that if $e(G)< \Lambda h$  for every plane-saturated subgraph $H$ of $G$, then Theorem~1.1 is immediate.	Thus we may suppose, for a contradiction, that there exists a plane-saturated subgraph $H$ with
	\begin{equation}\label{eq:contradiction}
		e(G)\ge \Lambda h.
	\end{equation}
	Since \(G\) has an edge,  $H$ is plane-saturated that implies \(h>0\). Throughout Sections~2 and~3, we work under the assumption \eqref{eq:contradiction}. We first explore the relationship between faces and $I$ of $H$.
	
\begin{lemma}\label{lem:universal-isolate}
	There are two vertices of $I$ lying in a  face of \(H\). Moreover, for
	every embedding \(\phi:H\to G\), the image set of $I$ is independent in \(G\).
\end{lemma}

\begin{proof}
	Let \(c\) and \(f\) be the numbers of components and faces of the embedded
	skeleton \(S\). Recall that \(S\) is obtained from \(H\) only by deleting isolated
	vertices. Hence, \(e(S)=h\). Euler's formula gives $
	s-h+f=1+c$. Hence, \(f=h-s+c+1\). 
	 Observe that every component of \(S\) contains an edge, so \(c\le h\).
	
	Suppose that every face contains at most one vertex of \(I\). Then, we have
	\[
	|I|\le f \ \text{and}\ 
	n=s+|I|\le s+f=h+c+1\le2h+1.
	\]
	By \eqref{eq:contradiction} and the value of $\Lambda$, we have that \(e(G)\ge\Lambda h\ge9\). Since \(G\) is
	simple, \(n\ge5\). Together with the planarity of $G$, this gives
	\[
	e(G)\le3n-6\le6h-3<9h\le\Lambda h,
	\]
	contrary to \eqref{eq:contradiction}. Thus there exists  some face of $H$ that contains two distinct
	vertices of I, say $u,v$.
	%\red{Under the assumption that every face contains at most one vertex of \(I\), we have \(n\le2h+1\). Since \(G\) is simple and \(e(G)\ge9\), we have \(n\ge5\).}
	% 修改原因：原句未明确说明 \(n\le2h+1\) 是在反设下得到的，且未说明 \(n\ge5\) 的由来。
	
	Let \(\phi:H\to G\) be an embedding. Assume to the contrary that \(\phi(I)\) contains an
	edge. Then, since every vertex of \(I\) is isolated in \(H\), we may permute
	the images of the vertices of \(I\), while fixing the images of
	\(V(H)\setminus I\), to obtain an embedding \(\phi':H\to G\)
	such that \(\phi'(u)\phi'(v)\in E(G)\). Thus, the edge \(rs\) can be drawn inside
	their common face of $H$, and \(\phi'\) then extends to an embedding
	\(H+uv\to G\), contradicting the plane saturation of $H$ with respect to $G$.
\end{proof}	

We now introduce a closure process on the matching components of \(S\).
	The process starts from the isolated vertices of \(H\) and repeatedly
	activates matching components that satisfy an eligibility condition defined
	below. It produces a terminal set \(Q\subseteq V(H)\), which will be used
	throughout the rest of the proof.
% 修改原因：在闭包过程正式定义之前，先交代目的、起点和产物 Q，给读者心理预期。

The closure process proceeds on $j\ge 0$. For $j+0$, set \(Q_0=I\).  Suppose that $Q_j$ has been defined for some $j\ge 0$. 
For $j+1$, if an unused matching component $ab$
satisfies
\begin{equation}\label{eq:eligible}
|N_{Q_j}(a)|=0	\ \text{and} \ |N_{Q_j}(b)|\ge2,
\end{equation}
then activate that component and put $Q_{j+1}=Q_j\cup\{a\}$.
Call $a$ its \emph{portal} and $b$ its \emph{mate}.

%Continue until no such ordered pair $a,b$ remains, and denotethe terminal set by $Q$.  

The process terminates as soon as no unused matching component satifying (\ref{eq:eligible}). We call it the \emph{closure process} and denote its terminal set
	by \(Q\). We call all vertices of $Q$ portals, including the initial vertices in $I$.

	\begin{lemma}\label{lem:pairwise-release}
		For every \(j\ge0\) and every two distinct vertices \(x,y\in Q_j\), there
		is an embedding $\phi:H\to G$ fixing every vertex outside
		$Q_j$ and mapping two distinct vertices of $I$ to $x,y$.
		Consequently, every $Q_j$, and in particular $Q$, is independent in $G$.
	\end{lemma}

	\begin{proof}
		We prove the assertion by induction on \(j\). Together with identity maping and Lemma \ref{lem:universal-isolate},  the assertion is immediate for \(Q_0=I\). % Note that the assertion is immediate for \(Q_0=I\): every vertex of \(I\) is isolated in \(H\), it thus has a trivial release	chain consisting of itself, which equals its support, and then any two  distinct vertices of \(I\) have	support-disjoint trivial chains.
		Assume the assertion holds for all \(Q_i\) with \(i\le j\). Let $ab$ be the unused matching component such that $Q_{j+1}=Q_j\cup \{a\}$. Let
		\(x,y\in Q_{j+1}\) be two distinct vertices.    If  \(x,y\in Q_j\), then, by the induction hypothesis, it is true. Since $Q_{j+1}\setminus Q_{j}=\{a\}$,  the remaining case is $x=a$ and $y\in Q_j$.  Let $z\in N_{Q_j(b)}\setminus\{b\}$. Note that $z,y\in Q_j$, hence,  there exists an embedding $\phi:H\to G$ fixing every vertex outside
		$Q_j$ and mapping two distinct vertices $u$ and $v$ of $I$ to $z,y$. Hence, $\phi(u)=z$ and $\phi(v)=y$, and then $\phi(a)=a$ and $\phi(b)=b$ by $a,b\notin Q_j$.  By exchange the images of $a$ and $u$, we	define \(\phi':H\to G\) by
		\[
		\begin{array}{c|ccc}
			w	& a & u & w\neq a,u \\ \hline
			\phi' & z & a& \phi(w)  
		\end{array}
		.\]  
		This is a bijection.  As $u$ is isolated and $ab$ is a matching
		component, the only changed edge image is $ab\mapsto zb\in E(G)$. 
		Thus $\psi'$ is the maping as required  for $a,y$, fixing every vertex
		outside $Q_{j+1}$.  
		In addition, if $G[Q_{j+1}]$ contains edges, then this contradicts
		Lemma \ref{lem:universal-isolate}.
		\end{proof}
		% 修改原因：用 x,y 作为待证对象贯穿证明，避免引入 a,a' 造成记号脱节；归纳假设改为对所有 Q_i（i≤j）成立；补出 L_{j+1}∩Q_j=∅ 说明三种 Case 互不重叠。
		
%		\begin{remark}\label{rem:pairwise}
%			Lemma~\ref{lem:pairwise-release} is a pairwise statement. For each pair
%			\(x,y\), the two release chains may be chosen anew; no globally fixed family
%			of pairwise support-disjoint release chains is asserted.
%		\end{remark}
%		
		
The following two lemmas are consequences of
Lemma~\ref{lem:pairwise-release}, which describe restrictions on the tree
components of \(S\).

\begin{lemma}\label{lem:q-leaf}
	Let \(T\) be a tree component of \(S\) with \(V(T)\cap Q=\varnothing\).
	Let \(a\) be a leaf of \(T\), and let \(b\) be its unique neighbor in
	\(T\). Suppose that \(N_Q(a)\) and \(N_Q(b)\)
	are nonempty. Then $
	N_Q(a)=N_Q(b)=\{q\}$ for some \(q\in Q\).
\end{lemma}
	
	\begin{proof}
		Suppose, for a contradiction, that there exist distinct
		\(q_1\in N_Q(a)\) and \(q_2\in N_Q(b)\). 
		Lemma~\ref{lem:pairwise-release} gives an embedding $\phi_1$ fixing every vertex outside $Q$
		 and mapping distinct $r_1,r_2\in I$ to $q_1,q_2$, respectively. Clearly, $\phi_1(r_1)=q_1$ and $\phi_2(r_2)=q_2$, and then $\phi_1(a)=a$ and $\phi_2(b)=b$ by $a,b\notin Q$.
		By exchange the images of $a$ and $r_2$, we define \(\phi'_1:H\to G\) by
			\[
		\begin{array}{c|ccc}
			w	& a & r_2 & w\neq a,r_2 \\ \hline
			\phi'_1 & z & a& \phi_1(w)  
		\end{array}
		.\]  
		This is a bijection. Observe that 
		Only the edge image of $ab$ 
		changes, becoming $q_2b\in E(G)$.  Hence, $\phi'_1$ is an embedding.  But it  maps
		the isolates $r_1,r_2$ to the adjacent vertices $q_1,a$,  contradicting
		Lemma~\ref{lem:universal-isolate}.
	\end{proof}
	\begin{lemma}\label{lem:p3-leaves}
		Let \(xyz\) be a \(P_3\)-component of \(S\). If
		\(q_1\in N_Q(x)\) and \(q_2\in N_Q(z)\), then \(q_1=q_2\).
	\end{lemma}

	\begin{proof}
		Suppose \(q_1\ne q_2\). By Lemma \ref{lem:pairwise-release}, we can obtain an embedding $\phi_2$ such that $\phi_2(r_1)=q_1$ and $\phi_2(r_2)=q_2$ for distinct $r_1,r_2\in I$. Observe that this component is vertex-disjoint from $Q$. We now define a new maping \(\phi'_2:V(H)\to V(G)\) by
			\[
		\begin{array}{c|ccccc}
			w	    & x & y &z & r_1 &w\neq x,y,z,r_1\\ \hline
			\phi'_2 & q_1 & x&y &z &\phi_2(w)  
		\end{array}
		.\]  
	 This is a cyclic
		permutation on \(\{u,v,w,q_1\}\), hence a bijection. Observe that, from $\phi_2$ to $\phi'_2$, the only change images of edges in $H$ are the two edges of the $P_3$-component.  
		\(H\) in this set are \(uv\) and \(vw\). In fact, their images are $
		q_1x, xy,$ under $\phi'_2$
		which lie in \(E(G)\). Thus \(\phi\) is an embedding.
		
		Note that  the isolated vertices \(q_1,q_2\in I\) are mapped to \(z,q_2\),
		respectively. Since \(q_2\in N_I(w)\), the images \(\phi(q_1)\) and
		\(\phi(q_2)\) are adjacent in \(G\), contradicting
		Lemma~\ref{lem:universal-isolate}. Hence \(q_1=q_2\).
	\end{proof}
	% 修改原因：用 uvw 与 Lemma 2.6 保持一致；用两行表格表示映射；补出映射是双射、边保持的理由；明确矛盾来自 q_2∈N_I(w) 保证 wq_2∈E(G)；结论直接写 q_1=q_2。

	\section{Counting lemmas}
	
	The structural restrictions of Section~2 are now used to estimate the
		number of vertices of \(S\) outside \(Q\). We use two ingredients: a
		weight on the vertices outside \(Q\) that measures how few neighbors they
		have in \(Q\), combined with the standard planar edge bounds; and the multiplicity
		parameters \(k_1,k_2\) to bound the degree-one and degree-two vertices.
	% 修改原因：承接第 2 节；明确目标；说明两套工具；不出现尚未定义的 F。

We first introduce some notation and defintions. %\red{Let $U\subsetneq V(G)$. For \(z\in V(G)\setminus Q\), define $\omega_U(z)=(2-d_U(z))_+.$}
Let \(ab\) be an unused matching component of \(S\). If both \(N_Q(a)\) and \(N_Q(b)\) are nonempty, then Lemma~\ref{lem:q-leaf} yields $
N_Q(a)=N_Q(b)=\{q\}
$ for a unique \(q\in Q\). In this case, \(ab\) is called \emph{tight}, and \(q\) is called its \emph{anchor}. Otherwise, by symmetry, assume that \(N_Q(a)=\varnothing\). Since the closure process has terminated, \(|N_Q(b)|\le 1\)~(Otherwise, \(ab\) would still be eligible for activation.). Thus every unused matching component \(ab\) satisfies $
|N_Q(a)|\le 1
\ \text{and}\ 
|N_Q(b)|\le 1.
$ %%可以不写

Let $t_1$ be the number of tight matching components, let \(X\subseteq Q\)
be the set of all distinct anchors, and put \(t_X=|X|\). Denoted by $U$ the endpoints of all tight matching components.  Then \(t_X\le t_1\).  Let \(t_2=|Q\setminus I|\), the number of activated matching components. Equivalently, \(|Q|=|I|+t_2\). %t量否与其它记号重合需要检查
Let \(t_3\) be the number of tree components of \(S\) which are neither
activated matching components nor tight matching components. We call them
\emph{residual tree components}. Observe that every residual tree component \(T\) is disjoint from
	\(Q\cup U\), and hence \(V(T)\subseteq R\).  %这句话可以不要，要的话要检查符号。

 Put
$
Q'=Q\setminus X, R=V(G)\setminus (Q\cup U)$ and $ r=|R|. %今天写到这
$
For \(i\in\{1,2\}\), set $
D_i=\{q\in Q': d_G(q)=i\}, d_i=|D_i|.
$
Finally, put $
P=Q'\setminus(D_1\cup D_2)
%F=G'-(D_1\cup D_2),$ and  $n_F=v(F)
.$ Since \(G\) has no isolated vertices, we have \(d_G(q)\ge3\) for each $q\in P$. 
For a center \(c\) of a residual star \(T\cong K_{1,\ell}\) with
\(\ell\ge3\) and \(d_P(c)=0\), let
$
t(c)=\bigl|\{y\in N_T(c):d_P(y)\ge2\}\bigr|.
$
Let \(\mathcal{D}\) be the set of all such centers \(c\) for which
\(t(c)\ge 3\).  
Based on quantities above, we now bound $|P|$ by the number of the residual tree components of $S$ and the order of $R$.

	\begin{lemma}\label{lem:weighted-planar}
	$
	|P|\le2r-t_3.
	$
	\end{lemma}
	% 修改原因：明确列出引理涉及的符号；用两句替代“Moreover”；保留“neighbor in G”的说明。

	\begin{proof}

		We first define a graph $B$ with vertex set $P\cup R$ whose edge set consists of all edges of $G$ between $P$ and $R$, together with, for each $c\in\mathcal D$, all edges $cy$ of the residual star centered at $c$ such that $y$ is a leaf with $d_P(y)\ge2$. 
		
		Observe that $E(B)\subseteq E(G)$. Hence, 
		\(B\) is a simple planar subgraph of
		\(G\). Since each vertex $c$
		of \(\mathcal{D}\) is \(P\)-zero, the added edges of \(G\) between \(P\) and \(R\)
		have their endpoints in \(P\) and \(R\setminus \mathcal{D}\), respectively; since the added edges of \(H\) are the form $cy$
		such that $c\in \mathcal{D}$ and $y\in R$ with \(d_P(y)\ge2\), those $y$ are not in \(\mathcal{D}\).
		Hence \(B\) is bipartite.  Thus $B$ is  bipartite, with bipartition
		$
		\bigl(P\cup\mathcal D,R\setminus\mathcal D\bigr).
		$

		Let \(Z_i\) be the set of vertices in the part \(R\setminus \mathcal{D}\) of
		\(B\) having degree \(i\) in \(B\), for \(i=0,1\). Since every added edge
		\(cy\) with \(c\in A_0\) satisfies \(y\in R\setminus \mathcal{D}\) and
		\(d_P(y)\ge2\), no such edge is incident with a vertex \(z\in
		R\setminus \mathcal{D}\) having \(d_P(z)\le1\). Hence the degree of \(z\) in
		\(B\) equals \(d_P(z)\), and therefore $
		Z_i=\{z\in R\setminus \mathcal{D}:d_P(z)=i\},$ for $ i=0,1.
		$
	
		By the choice of $P$ and $\mathcal{D}$, we derive that $d_B(z)\ge 3$ for every $z\in P\cup \mathcal{D}$. Deleting $D_1\cup D_2$ from $B$ removes exactly $|Z_1|$ edges. Deleting $Z_0\cup Z_1$ removes exactly $|Z_1|$ edges.
		Recall that since \(B\) is a simple bipartite planar graph, we may use the coarse bound $
		e(B)\le 2v(B)$ from Proposition~\ref{prop:planar-bounds}. 
		 Hence, this leads to
		\begin{equation*}%\label{eq:B-lower}
			3|P|+\sum_{c\in\mathcal D}t(c)
			\le e(B)\le2(|P|+r)-2|Z_0|-|Z_1|.
		\end{equation*}
	Equivalently, 
		\begin{equation}\label{eq:B-lower}
			|P|\le2r-\left(2|Z_0|+|Z_1|+\sum_{c\in\mathcal D}t(c)\right).
		\end{equation}

		With \(Z_0,Z_1\), and \(\mathcal D\) obtained from the construction of \(B\), define
		\(\chi:R\to\mathbb R_{\ge0}\) by
		\[
		\chi(z)=
		\begin{cases}
			2, & z\in Z_0,\\
			1, & z\in Z_1,\\
			t(z), & z\in\mathcal D,\\
			0, & \text{otherwise}.
		\end{cases}
		\]
		Since \(Z_0,Z_1\), and \(\mathcal D\) are pairwise disjoint, the preceding identity gives
		\[
		W=2|Z_0|+|Z_1|+\sum_{c\in\mathcal D}t(c)
		=\sum_{z\in R}\chi(z).
		\]
		For a residual tree component \(T\), put
		$
		\chi(T)=\sum_{z\in V(T)}\chi(z).
		$
		Each residual tree component lies in \(R\), and distinct residual tree
		components are vertex-disjoint. Since \(\chi(z)\ge0\) for every \(z\in R\),
		we have
		$
		W\ge \sum_T \chi(T),
		$
		where the sum is over all residual tree components. Thus it is enough to
		prove that \(\chi(T)\ge3\) for every residual tree component \(T\).
		
		{\bf Clam 1}  \(\chi(T)\ge3\).

	{\bf Case 1:} \(T\) is not a star. Then the diameter of $T$ is at least three, which implies that \(T\) has two leaves whose
	neighbors in \(T\) are distinct. Hence, $T$ contains two disjoint leaf-neighbor pairs. For each pair \(\{a,b\}\) of $T$, where
	\(a\) is a leaf and \(b\) its neighbor, either both \(a\) and \(b\) have
	neighbors in \(P\), in which case Lemma~\ref{lem:q-leaf} gives
	\(N_P(a)=N_P(b)=\{q\}\) for some \(q\in P\), so \(d_P(a)=d_P(b)=1\) and
	\(\chi(a)+\chi(b)=2\); or at least one of them has no neighbors in $P$,
	so \(\chi(a)+\chi(b)\ge2\). Thus \(\chi(T)\ge4\).
	
	{\bf Case 2:} \(T\cong K_{1,\ell}\) with \(\ell\ge3\). Let \(c\) be the
	center, and let \(\ell_i\) be the number of leaves \(z\) with \(d_P(z)=i\)
	for \(i=0,1\), while \(\ell_2\)($=t(c)$) is the number of leaves with \(d_P(z)\ge2\). In particular,  $c\in \mathcal{D}$ if $\ell_2\ge 3$.
	
	We first assum that \(d_P(c)=0\). If $\ell_2\ge3$,  then it is trivial. So we assume that $\ell_2\le 2$. This implies that $c\in D_0$, it follows that $\chi(T)\ge 2+1=3$.

	If \(d_P(c)=1\), then $c\in D_1$, and Lemma~\ref{lem:q-leaf} brings that  each leaf $z$ of $T$ has $d_P(z)\le 1$. Then $
	\chi(T)\ge 4.
	$
	
	If \(d_P(c)\ge2\), then Lemma~\ref{lem:q-leaf} implies that every leaf is $P$-zero
	, so \(\chi(T)\ge2\ell>3\).
	
	{\bf Case 3:} \(T\cong P_3\), say \(T=abc\).  If at most one vertex of $T$ is $P$-positive, then $\chi(T)\ge 3$. Otherwise, there are at least two $P$-positive vertices. by Lemma~\ref{lem:q-leaf}, they are belonging to $D_1$. This implies that $\chi(T)\ge 3$.
		
	{\bf Case 4:} \(T\cong K_2\). Let $T=ab$. Since \(T\) is a residual tree component of $S$, it is an unused matching component that is not tight. By Lemma~\ref{lem:q-leaf},
 we can asssume that $d_Q(a)=0$, then
 \eqref{eq:eligible} gives \(d_Q(b)\le1\), and therefore
 $\chi(T)\ge 3$. 
%Let $T=ab$. Since \(T\) is a residual tree component of $T$, it is an unused matching component that is not tight, at most one of $a, b$ has a neighbor in $Q$. Otherwise,  By Lemma~\ref{lem:q-leaf},   $a$ and $b$
%have a unique common neighbor in $Q$. This implies that $T$ is tight, a contradiction.	
%Asssume that $d_Q(a)=0$, then
%\eqref{eq:fixedpoint} gives \(d_Q(b)\le1\), and therefore
%$	\omega_Q(a)+\omega_Q(b)\ge3.$ Together with \(P\subseteq Q\),
%this 	gives \(W(T)\ge3\).  
	
Combining cases above, the proof is complete.	
\end{proof}

	We next estimate the possible neighborhoods of vertices in \(D_1\) and
	\(D_2\). For \(D_1\), these neighborhoods are single vertices of \(R\);
	for \(D_2\), they are edges of the auxiliary neighborhood graph on \(R\).
	In both cases the residual tree components reduce the effective count by
	\(t_3\), and the neighborhood multiplicity assumptions yield the two
	bounds below.

%	\red{We next estimate the possible neighborhoods of vertices in \(D_1\)
%		and \(D_2\). For \(D_1\), these neighborhoods are single vertices of
%		\(R'\); for \(D_2\), they are two-element subsets of \(R'\). In both cases
%		the residual tree components contribute the factor \(|R'|-\tau\), and the
%		multiplicity parameters \(k_1\) and \(k_2\) yield the two bounds below.}
	% 修改原因：只说明 \(D_1,D_2\) 的邻域类型、residual tree components 贡献的因子、多重度参数的作用；不解释符号定义。

%	\red{We next estimate the possible neighborhoods of vertices in \(D_1\)
%		and \(D_2\). For \(D_1\), these neighborhoods are single vertices of
%		\(R'\); for \(D_2\), they are two-element subsets of \(R'\). The residual
%		tree components contribute the factor \(|R'|-\tau\), and the multiplicity
%		parameters \(k_1\) and \(k_2\) yield the two bounds below.}
%这一段后面再修改

\begin{lemma}\label{lem:D1} $
	|D_1|\le k_1(r-t_3)
	$ and $
	|D_2|\le 3k_2(r-t_3).
	$
\end{lemma}
% 修改原因：明确列出引理涉及的符号；去掉“We have”。

\begin{proof}
Recall that \(D_i=\{q\in Q':d_G(q)=i\}\) for $i=1,2$. Moreover, for each $q\in D_1\cup D_2$, $N_G(q)\subseteq R$. We now show the first inequality.  Let $
A_1=N_G(D_1).
$
Then every $q\in D_1$ satisfies $N_G(q)=\{z\}$ for some $z\in A_1$.  For each \(z\in A_1\), at most \(k_1\)
vertices of degree one in \(G\) have neighborhood \(\{z\}\) by the
definition of \(\mathcal G_{k_1,k_2}\). Hence, we obtain that 
 $|D_1|\le k_1|A_1|$.
% 修改原因：先说明 \(D_1\) 的定义和度为 1，再用 Lemma 3.1 说明邻居在 \(R'\) 中，再定义 \(A_1\)，最后用多重度参数。
	
	We  claim that every residual tree component of $S$ contains a vertex outside
	\(A_1\). Suppose, to the contrary, that there exists 	residual tree component \(T\)  such that \(V(T)\subseteq A_1\). Let \(a\) be a leaf of \(T\), and let \(b\)
	be its neighbor in \(T\). Since \(a,b\in A_1=N_G(D_1)\), there exist
	\(q_a,q_b\in D_1\) with \(q_a a\in E(G)\) and \(q_b b\in E(G)\). On the other hand, by
	Lemma~\ref{lem:q-leaf}, we derive that \(N_Q(a)=N_Q(b)=\{q\}\) with $q_a=q_b=q$. Hence \(d_G(q)\ge2\). But this contradicts the condition \(q\in D_1\). Therefore \(V(T)\setminus A_1\ne\varnothing\)
	for every residual tree component \(T\).
	
	Since the residual tree components are vertex-disjoint,   \(R\setminus A_1\) contains at least \(t_3\)
	distinct vertices from all residual tree components. Hence \(|R\setminus A_1|\ge t_3\),
	i.e.\ \(|A_1|\le|R|-t_3\). Combining this with \(|D_1|\le k_1|A_1|\)
	yields \(|D_1|\le k_1(|R|-t_3)\).

We now turn to the second estimate. Recall that every $q\in D_2$ satisfies $d_G(q)=2$ and $N_G(q)\subseteq R$. Thus $N_G(q)$ is a two-element subset of $R$. Set $\mathcal N_2=\{N_G(q):q\in D_2\}$. For each $X=\{x,y\}\in\mathcal N_2$, fix a representative $q_X\in D_2$ such that $N_G(q_X)=X$. Define an auxiliary graph $O$ with vertex set $R$ by putting $xy\in E(O)$ if and only if $\{x,y\}\in\mathcal N_2$. For each edge $xy\in E(O)$, the path $xq_Xy$, where $X=\{x,y\}$, is the corresponding $P_3$ in $G$. Let $G_{P_3}$ be the subgraph of $G$ consisting of all these paths, together with the remaining vertices of $R$ as isolated vertices. Then $G_{P_3}$ is obtained from $O$ by subdividing each edge exactly once. The internal vertices $q_X$ are pairwise distinct and lie outside $R$, while different paths may share only endpoints in $R$. Suppressing all internal vertices $q_X$ therefore recovers $O$. Since $G_{P_3}$ is planar, $O$ is planar.
	
For each \(X\in\mathcal N_2\), at most \(k_2\) vertices of degree two in
\(G\) have neighborhood \(X\), by the definition of
\(\mathcal G_{k_1,k_2}\). Hence, 
	\begin{equation}\label{eq:D2-multiplicity}
		|D_2|\le k_2|\mathcal N_2|=k_2e(O).
	\end{equation}

	We show that every residual tree component $T$ contains a component of $O$ on one or two vertices. If $T$ has a vertex which is isolated in $O$, then this vertex is such a component. Otherwise every vertex of $T$ is incident with an edge of $O$. Choose a leaf $a$ of $T$, and let $b$ be its unique neighbor in $T$. Then both $a$ and $b$ are incident with edges of $O$, and hence both have neighbors in $D_2\subseteq Q$. By Lemma~\ref{lem:q-leaf}, $N_Q(a)=N_Q(b)=\{q\}$ for some $q\in Q$. Since an edge of $O$ incident with $a$ comes from a vertex of $D_2$ adjacent to $a$, this vertex must be $q$. Thus $q\in D_2$. As $d_G(q)=2$ and $q$ is adjacent to both $a$ and $b$, we have $N_G(q)=\{a,b\}$. Therefore $ab\in E(O)$. Moreover, every edge of $O$ incident with $a$ or $b$ comes from a vertex of $D_2$ in $N_Q(a)$ or $N_Q(b)$, respectively. Since both neighborhoods are equal to $\{q\}$, the edge $ab$ is an isolated $K_2$-component of $O$.
	
	Choose one such component of $O$ from each residual tree component. These chosen components are distinct because the residual tree components are vertex-disjoint. Hence $c(O)\ge t_3$. By \eqref{eq:component-planar},
	\[
	e(O)\le 3\bigl(v(O)-c(O)\bigr)\le 3(r-t_3).
	\]
	Together with \eqref{eq:D2-multiplicity},%% 需要用到平面图的边数估计，这是需要在第二章介绍还是直接在这引用
	 this proves the second bound.	
\end{proof}

	\section{Proof of Theorem~\ref{thm:main}}
	
	\begin{proof}	
	Let $H$ be an arbitrary spanning plane-saturated subgraph of $G$. Suppose that \eqref{eq:contradiction} holds, and apply the preceding constructions.
	
	We first estimate $r=|R|$ by considering the components of the skeleton $S$. An activated matching component contributes only its mate to $R$. A tight matching component contributes no vertex to $R$, since both of its endpoints are deleted. A residual tree $T$ contributes $v(T)=e(T)+1$, while every component of $S$ containing a cycle contributes at most as many vertices as edges. Hence
	\begin{equation}\label{eq:r-count}
		0\le r-t_3\le h-t_2 .
	\end{equation}
	The lower bound follows because the $t_3$ residual trees are vertex-disjoint nonempty subgraphs contained in $R$. This count also covers the case in which an activated portal is an anchor: the portal lies in $Q$, while its mate remains in $R$.
	
	Put $d_i=|D_i|$ for $i=1,2$, and let $U$ be the set of endpoints of the tight matching components. Thus $|U|=2t_2$. Thus \(V(G)=P\cup D_1\cup D_2\cup X\cup R\cup U\), and this union is disjoint.  Since $D_1\cup D_2$ is independent, deleting it removes exactly $d_1+2d_2$ edges. The remaining planar graph has $|P|+r+x+2t_2$ vertices and is nonempty, because $h>0$ and the two endpoints of an edge of $H$ cannot both lie in the independent set $Q$. Therefore, by \eqref{eq:component-planar}, %% 需要用到平面图的边数估计，这是需要在第二章介绍还是直接在这引用
	\[
	e(G)\le 3(|P|+r+x+2t_2)+d_1+2d_2-3 .
	\]
	
	By Lemmas~\ref{lem:weighted-planar} and ~\ref{lem:D1}, and then by \eqref{eq:r-count} and $x\le t_2$, we obtain
	\begin{align*}
		e(G)
		&\le 9(r-t_3)+d_1+2d_2+3x+6t_2-3\\
		&\le \Lambda(r-t_3)+3x+6t_2-3\\
		&\le \Lambda(h-t_2)+9t_2-3\\
		&=\Lambda h-(\Lambda-9)t_2-3\\
		&<\Lambda h.
	\end{align*}
	This contradicts \eqref{eq:contradiction}. Thus every spanning plane-saturated subgraph $H$ satisfies
	\[
	\frac{e(H)}{e(G)}>\frac1\Lambda .
	\]
	By Lemma~\ref{lem:reductions}, the same conclusion holds for all plane-saturated drawings and for the original host graph with its isolated vertices. This proves the theorem.	
	\end{proof}

	\noindent{\bf Acknowledgments}
		
		Qing Cui is partially supported by the National Natural Science Foundation of China (No. 12171239). Shengjin Ji is partially supported by the National Natural Science Foundation of China (No. 12561061).

\end{document}